\documentclass[a4paper, 11pt]{amsart}
\usepackage{mathrsfs}
\usepackage{amsmath}
\usepackage{amssymb}
\usepackage{verbatim}
\usepackage{xcolor}
\usepackage{CJK}
\usepackage{graphicx}
\usepackage[all]{xy}
\usepackage{appendix}
\DeclareGraphicsExtensions{.eps,.ps,.jpg,.bmp}

\theoremstyle{plain}

\newtheorem{theorem}{Theorem}[section]
\newtheorem{proposition}{Proposition}[section]

\theoremstyle{definition}

\theoremstyle{remark}

\newcommand{\eps}{\epsilon}

\newcommand{\Tr}{\operatorname{Tr}}
\newcommand{\Id}{\operatorname{Id}}

\title{Limiting behavior of a class of Hermitian Yang-Mills metrics, II:
exponential decay}
\author{Jixiang Fu}

\begin{document}

\maketitle

\begin{abstract}
In the note, the geometric set-up, the rank two bundle, the local HYM ansatz, and the
global gluing construction are the same as in the preceding work \cite{Fu}.
The new point is an exponential estimate for the radial ordinary
differential equation obtained near each branch point.  If
$u_\eps$ denotes the local radial solution and
$\frac12\ln r$ the singular limiting solution, then for every
integer $k\ge0$, there exist positive constants $C_k$ and $c_k$ such that
\[
\big\|
u_\eps-\frac12\ln r
\big\|_{C^k([r_0,2r_0])}
\le C_k e^{-c_k/\eps}.
\]
Consequently, all results of the preceding paper can be refined.
\end{abstract}

\section{Modified statements of the main results}

We keep the notation of the preceding work \cite{Fu}.

The K\"ahler manifold $X$ is the product $B\times T$
of two copies of the complex one-torus  $\mathbb C/\Gamma$, where
$\Gamma=\mathbb Z+i\mathbb Z$, and the family of product
metrics $\omega_\epsilon$  are flat and  have areas $\epsilon$
and $\epsilon^{-1}$ on  $T$ and $B$ respectively.

The holomorphic vector bundle $V$ over $X$ is
constructed as follows (cf. \cite{Fr,FMW}). Let $T^\ast$ be the dual
of $T$ and let $X^\ast=T^\ast\times B$.  Let  $Y$ be a compact
(complex) curve of $X^\ast$ such that the induced projection
$\varphi:Y\rightarrow B$ is a two-sheet branched cover with $n$ branched points. Denote the other induced map by $q: Y\to T^\ast$.
Denote
\begin{equation*}\label{171221}
\iota=(q,\textup{id}_{T}):Y\times T \to T^\ast\times T,
\qquad p_2=(\varphi,\textup{id}_T):Y\times T\to X
\end{equation*}
and denote by $p_1$ the projection map from $Y\times T$ to $Y$.
Let $\mathcal P$ be the Poincar\'e
line bundle on $T^\ast\times T$.
 Then for any
degree zero line bundle $\mathcal F$ over $Y$, we can form a line bundle over $Y$
\begin{equation*}
\mathcal N=K^{1/2}_Y\otimes \varphi^\ast K_B^{-1/2}\otimes
\mathcal F
\end{equation*}
and a rank two vector bundle over $ X$ with zero degree
\begin{equation*}
V=p_{2\ast}(\iota^*\mathcal P\otimes p_1^*\mathcal N).
\end{equation*}
By an adiabatic argument (cf. \cite{FMW}), $V$ is $\omega_\epsilon$-slope stable for
small $\epsilon$. Hence, by the Donaldson-Uhlenbeck-Yau theorem,  there exists  a family of irreducible
Hermitian Yang-Mills (HYM) metrics $H_{1,\epsilon}$ on $V$ with respect to
$\omega_\epsilon$. As a consequence of $c_1(V)=0$,
  the associated curvature forms $\Theta(H_{1,\epsilon})$ satisfy
\begin{equation*}\label{06}
\Lambda_{\omega_\epsilon}\Theta(H_{1,\epsilon})=0.
\end{equation*}
The approximate metric constructed in \cite{Fu} is
denoted by $H_{0,\eps}$, and
\[
H_{1,\eps}(\cdot,\cdot)=H_{0,\eps}(H_\eps\cdot,\cdot),
\qquad
\det H_\eps=1.
\]

The new point is an exponential estimate for the radial ordinary
differential equation obtained near each branch point.  If
$u_\eps$ denotes the local radial solution and
$\frac12\ln r$ the singular limiting solution of equation (\ref{00}), then for every
integer $k\ge0$, there exist positive constants $C_k$ and $c_k$ such that
\[
\big\|
u_\eps-\frac12\ln r
\big\|_{C^k([r_0,2r_0])}
\le C_k e^{-c_k/\eps}.
\]
Consequently,
The following theorem replaces Theorem 2 of the preceding work \cite{Fu}.

\begin{theorem}[Exponentially accurate approximate metric]
For every non-negative integer $k$, there exist positive constants
$C_k$ and $c_k$ such that, for all sufficiently small $\eps>0$,
\[
\left\|
\Lambda_{\omega_\eps}\Theta(H_{0,\eps})
\right\|_{C^k}
\le C_k e^{-c_k/\eps}.
\]
Moreover,
\[
\Tr\Lambda_{\omega_\eps}\Theta(H_{0,\eps})=0.
\]
\end{theorem}

The conditional higher order estimate becomes the following exponential
version of Theorem 3 in \cite{Fu}.

\begin{theorem}[Conditional higher order estimate]
Fix a non-negative integer $k$.  Assume that there exist positive constants
$C_0$ and $c_0$ such that for all sufficiently small $\eps$,
\[
\|H_\eps-\Id\|_{C^0}\le C_0e^{-c_0/\eps}.
\]
Then there exist positive constants $C_k$ and $c_k$ such that
\[
\|H_\eps-\Id\|_{C^k}\le C_k e^{-c_k/\eps}.
\]
\end{theorem}

The lower bound estimate is also improved.

\begin{theorem}[Exponential lower bound]
There exist positive constants $C$ and $c$ such that, for all sufficiently small
$\eps>0$,
\[
\inf_{x\in X}\Tr H_\eps(x)\le 2+Ce^{-c/\eps}.
\]
\end{theorem}

\section{Replacement for Section 4: the ODE estimate}

We consider the Dirichlet problem
\begin{equation}\label{00}
\begin{cases}
\Delta u_\eps
=
4\pi^2\eps^{-2}
\left(e^{2u_\eps}-r^2e^{-2u_\eps}\right)
&\text{in }B_{2r_0}(0),\\[4pt]
u_\eps=\frac12\ln(2r_0)
&\text{on }\partial B_{2r_0}(0).
\end{cases}
\end{equation}
The singular solution is
\[
u_s(r)=\frac12\ln r.
\]

As in the preceding work, after substituting $\overline u$ for $2u-\ln(2r_0)$, $x_1$ for $\frac{x_1}{2r_0}$,
$x_2$ for $\frac{x_2}{2r_0}$, $r^2$ for $\frac{r^2}{4r_0^2}$, and $\eps$ for $\frac{\eps}{8\pi}r_0^{- 3/ 2}$,
the problem is reduced to
\[
\begin{cases}
\Delta \overline u_\eps
=
\eps^{-2}
\bigl(e^{\overline u_\eps}-r^2e^{-\overline u_\eps}\bigr)
&\text{in }B_1(0),\\[4pt]
\overline u_\eps=0
&\text{on }\partial B_1(0).
\end{cases}
\]
By the same uniqueness and symmetry argument, $\overline u_\eps=\overline u_\eps(r)$ and
\[
\overline u_\eps''(r)+r^{-1} \overline u_\eps'(r)
=
\eps^{-2}
\bigl(e^{\overline u_\eps(r)}-r^2e^{-\overline u_\eps(r)}\bigr),
\qquad
\overline u_\eps(1)=0.
\]

Define
\[
v_\eps(r)=\overline u_\eps(r)-\ln r.
\]
Then
\begin{equation}\label{12}
v_\eps''+r^{-1} v_\eps'
=
2\eps^{-2}r\sinh v_\eps,
\qquad
v_\eps(1)=0.
\end{equation}
The known monotonicity properties are
\[
v_\eps>0,\qquad v_\eps'<0,\qquad v_\eps''>0
\quad\text{on }(0,1).
\]
We also recall the elementary integral estimate
\[
\sinh v_\eps(1/4)\le C\eps^2,
\]
hence, since $v_\eps>0$ and $v_\eps\le \sinh v_\eps$ for small
$v_\eps$,
\begin{equation}\label{14}
v_\eps(1/4)\le C\eps^2.
\end{equation}

Set
\[
q_\eps(r)
=
2\eps^{-2}r\frac{\sinh v_\eps(r)}{v_\eps(r)}.
\]
Then (\ref{12}) is
\[
v_\eps''(r)+r^{-1} v_\eps'(r)-q_\eps(r)v_\eps(r)=0.
\]
Since $\sinh t/t\ge1$ for $t>0$, on $[1/4,1]$ we have
\begin{equation}\label{13}
q_\eps(r)\ge 2\eps^{-2}r\ge \frac12\eps^{-2}.
\end{equation}

Let
\[
L_\eps w(r)=w''(r)+r^{-1} w'(r)-q_\eps(r)w(r).
\]
If $L_\eps w\ge0$ on $(a,1)$, $w(a)\le 0$ and $w(1)\le0$, then
$w\le0$ on $[a,1]$.  Indeed, a positive interior maximum would satisfy
$w'=0$, $w''\le0$, and hence
\[
L_\eps w\le -q_\eps w<0,
\]
a contradiction.

Take
\[
a=\frac14,\qquad \mu=\frac1{2\eps},
\]
and define
\[
\Phi_\eps(r)
=
v_\eps(a)
\frac{\sinh\bigl(\mu(1-r)\bigr)}
{\sinh\bigl(\mu(1-a)\bigr)}.
\]
Then
\[
\Phi_\eps(a)=v_\eps(a),\qquad
\Phi_\eps(1)=0,
\qquad
\Phi_\eps''=\mu^2\Phi_\eps.
\]
Moreover $\Phi_\eps'<0$, and by (\ref{13}),
\[
L_\eps\Phi_\eps
=
(\mu^2-q_\eps)\Phi_\eps+r^{-1}\Phi_\eps'
\le
\bigl(\frac1{4\eps^2}-\frac1{2\eps^2}\bigr)\Phi_\eps
<0.
\]
Therefore
\[
L_\eps(v_\eps-\Phi_\eps)>0,
\qquad
(v_\eps-\Phi_\eps)(a)=(v_\eps-\Phi_\eps)(1)=0.
\]
The comparison principle gives
\[
0\le v_\eps(r)\le \Phi_\eps(r),\qquad r\in[a,1].
\]
By direct computation,
we obtain
\[
\frac{\sinh\bigl(\mu(1-r)\bigr)}
{\sinh\bigl(\mu(1-a)\bigr)}
\le
C e^{-\mu(r-a)}.
\]
Thus
\begin{equation}\label{15}
v_\eps(r)
\le
C v_\eps(a)e^{-\frac{r-a}{2\eps}}.
\end{equation}
Combining (\ref{14}) with (\ref{15}), for $r\in[1/2,1]$,
\begin{equation}\label{16}
0\le v_\eps(r)
\le C\eps^2e^{-\frac1{8\eps}}
\le C e^{-\frac1{16\eps}}.
\end{equation}

For higher derivatives, fix $r_0\in[1/2,1]$ and put
\[
\rho=\frac{r-r_0}{\eps},
\qquad
V_\eps(\rho)=v_\eps(r_0+\eps\rho).
\]
Then
\[
V_\eps''+
\frac{\eps}{r_0+\eps\rho}V_\eps'
=
2(r_0+\eps\rho)\sinh V_\eps .
\]
On every fixed compact $\rho$-interval the coefficients are uniformly
bounded, and the $C^0$ norm of $V_\eps$ is exponentially small by
(\ref{16}).  Standard interior estimates for ordinary differential equations
give
\[
\left|
\frac{d^mV_\eps}{d\rho^m}(0)
\right|
\le
C_m\|V_\eps\|_{C^0}.
\]
Consequently,
\[
|v_\eps^{(m)}(r_0)|
\le
C_m\eps^{-m}e^{-c/\eps}
\le
C_m e^{-c_m/\eps}.
\]
Since $r_0\in[1/2,1]$ is arbitrary, we have proved:

\begin{proposition}[Exponential ODE estimate]
For every non-negative integer $k$ there exist constants $C_k,c_k>0$
such that
\[
\|\overline u_\eps-\ln r\|_{C^k([1/2,1])}
\le C_k e^{-c_k/\eps}.
\]
Returning to the original variables gives
\begin{equation}\label{11}
\big\|
u_\eps-\frac12\ln r
\big\|_{C^k([r_0,2r_0])}
\le
C_k e^{-c_k/\eps}.
\end{equation}
\end{proposition}

The interior polynomial estimates for $\|u_\eps\|_{C^k(B_R(0))}$ used
later in the proof are unchanged.

\section{Replacement for the end of Section 5}

We keep the construction of $\textbf h_\eps$ and the conformal normalization
\[
H_{0,\eps}=e^{-\frac12(\phi_1+\phi_2)}\textbf h_\eps
\]
unchanged.  The trace identity
\[
\Tr\Lambda_{\omega_\eps}\hat\Theta(H_{0,\eps})=0
\]
is also unchanged.

On the gluing annuli one has, as before,
\[
\frac{i}{2}\Lambda_{\omega_\eps}\hat\Theta( H_{0,\eps})
=
\psi_\eps
\begin{pmatrix}
1&0\\
0&-1
\end{pmatrix},
\]
where
\[
\psi_\eps
=
\frac{\pi^2}{\eps}r(\phi_\eps-\phi_\eps^{-1})
-\frac{\eps}{2}\partial_z\partial_{\bar z}\ln\phi_\eps,
\]
and
\[
\phi_\eps
=
\frac{
1+\rho(e^{-(u_\eps-\frac12\ln r)}-1)
}{
1+\rho(e^{u_\eps-\frac12\ln r}-1)
}.
\]
By (\ref{11}),
\[
\|\psi_\eps\|_{C^k}
\le
C_k\eps^{-1}e^{-c_k/\eps}
\le
C_k e^{-c'_k/\eps}.
\]
Therefore,
\[
\left\|
\Lambda_{\omega_\eps}\Theta(H_{0,\eps})
\right\|_{C^k}
\le
C_k e^{-c_k/\eps}.
\]

\section{Replacement for Section 6}

The normalization and the scalar inequality remain the same.  Namely,
after normalizing $H_{1,\eps}$ by $\det H_\eps=1$, one has
\[
-\Delta_{\omega_\eps}\ln\Tr H_\eps
\le
4\|\Lambda_{\omega_\eps}\Theta(H_{0,\eps})\|_{C^0}.
\]
By (5.2),
\[
-\Delta_{\omega_\eps}\ln\Tr H_\eps
\le
Ce^{-c/\eps}.
\]
The Moser iteration argument gives
\[
\sup_X\ln\Tr H_\eps
\le
(1+Ce^{-c/\eps})
\int_X\ln\Tr H_\eps\,\frac{\omega_\eps^2}{2!}.
\]
No upper $C^0$ estimate is concluded from this inequality alone.

\section{Replacement for Section 7: lower bound}

The proof of the lower bound is unchanged except that every occurrence
of the polynomial curvature error
\[
\|\psi_\eps\|_{C^0}\le C\eps^N
\]
is replaced by the exponential estimate
\[
\|\psi_\eps\|_{C^0}\le Ce^{-c/\eps}.
\]
All possible losses in the proof are polynomial powers of $\eps$, and
therefore are absorbed into the exponential term.  Consequently, one
obtains Theorem 3.

\section{Replacement for Section 8: conditional higher order estimate}
The higher order argument of the preceding work is also unchanged in
structure.  The only inputs are:

\[
\|H_\eps-\Id\|_{C^0}
\]
and the curvature error estimate for $H_{0,\eps}$.

Assume now that
\[
\|H_\eps-\Id\|_{C^0}\le C_0e^{-c_0/\eps}.
\]
In the notation of the preceding work, the equation for
$\mathcal H_\eps=H_\eps-\Id$ has the schematic form
\[
\Delta_\eps\mathcal H_\eps
=
Q(\nabla_\eps\mathcal H_\eps,\nabla_\eps\mathcal H_\eps)
+
P(\nabla_\eps\mathcal H_\eps)
+
R_\eps,
\]
where
\[
\|R_\eps\|_{C^j_\eps}\le C_j\eps^{-N_j}e^{-c/\eps}
\le C_j e^{-c_j/\eps}.
\]
The proof using the Gagliardo-Nirenberg inequality and the rescaled
coordinates gives, for every $m\ge0$,
\[
\|\nabla_\eps^m\mathcal H_\eps\|_{L^2}
\le
C_m e^{-c_m/\eps}.
\]
Sobolev embedding in the rescaled coordinates then yields
\[
\|\mathcal H_\eps\|_{C^k}
\le
C_k e^{-c_k/\eps}.
\]
Thus Theorem 2 holds.

\address{
\small Current address: School of Mathematical Sciences, Fudan University, Shanghai 200433, People's Republic of China.
\small E-mail address: majxfu@fudan.edu.cn,
}

\begin{thebibliography}{99}
\bibitem{Fr} R. Friedman, {\sl Rank two vector bundles over regular elliptic
surfaces}, Invent. Math. {\bf 96}(1989), 283-332.
\bibitem{FMW}R. Friedman, J. Morgan and E. Witten, {\sl Vector
bundles and F theory},  Comm. Math. Phys.  {\bf 187}(1997), 679-743.
\bibitem{Fu} J. Fu. Limiting behavior of a class of
Hermitian Yang-Mills metrics, I. Sci. China
Math. {\bf 62}(2019), 2155-2194.

\end{thebibliography}
\end{document}